\documentclass{article}
\usepackage{authblk}
\usepackage{amsthm}
\usepackage{amssymb}
\usepackage{amsmath}
\usepackage{amsfonts}
\usepackage{wrapfig}
\usepackage{rotating}
\usepackage{pgfpages}
\theoremstyle{plain}
\newtheorem{theorem}{Theorem}

\theoremstyle{definition}
\newtheorem{definition}[theorem]{Definition}

\begin{document}
\title{Dynamical systems' models for the prediction of multi-variable time series. Wikipedia's traffic example.}
\author[Victoria Rayskin]{Victoria Rayskin\thanks{The author is very grateful to A. Halfaker of Wikimedia Foundation for his help and advice with the data.\\
{Research was sponsored by the Army Research Office and was
accomplished under Grant Number W911NF-19-1-0399. The views and conclusions contained in this
document are those of the authors and should not be interpreted as representing the official policies, either
expressed or implied, of the Army Research Office or the U.S. Government. The U.S. Government is
authorized to reproduce and distribute reprints for Government purposes notwithstanding any copyright
notation herein.}} \\ victoria.rayskin@tufts.edu}

%\address{Department of Mathematics\\ Tufts University\\ Medford, MA 02155-5597}
%\email{victoria.rayskin@tufts.edu}
\date{\today}
\maketitle
\begin{abstract}
The models VAR, ARIMA, Holt-Winters, are frequently used for short-term forecasts of multivariate time series.
In this paper we consider models constructed with the help of dynamical systems that have relatively simple limiting behavior. Switching between different trajectories of the phase portrait, we obtain a high precision prediction. Moreover, the dynamical system approach provides the global qualitative picture of the model's phase portrait, and allows us to discuss multidimensional patterns and long-term properties of the process. The simple limiting behavior allows us to associate different trends with different process's realization scenarios that can be influenced by externalities. 

We demonstrate these ideas using the examples of the Wikipedia's traffic of Readers, Contributors and Edits. First, we consider the two-dimensional model, predicting the traffic of Readers and Edits. The prediction precision is higher than the two-dimensional VAR prediction. Different trends (corresponding to different fixed points) can be associated with different platform's incentives.
Then, adding the Contributors data, we discuss the three-dimensional model (more precise than the three-dimensional VAR). It provides a more accurate short-term prediction of Edits than the two-dimensional dynamic model. The global picture shows that the number of new Edits tends to decline in the future, while the number of new Contributors and Readers will grow in the long run.
\end{abstract}

\section{Introduction}\label{intro}
Modern availability of large multidimensional data sets creates opportunities for the development of new methods of data analysis. Traditional time series models (e.g., ARIMA, Vector Auto Regression, Holt-Winters) fit data into a single one- or multi-dimensional trajectory. Instead, we can try to fit the data into a vector field, which provides infinitely many trajectories for various initial states. We can move from one trajectory to another for a better fitting. In some cases, if traditional time series models are used, the complexity of the fitted trajectory can be significant. However, the same data can be accurately fitted into relatively simple differential equations, which are non-chaotic (according to R. Devaney definition of chaos). For example, the dynamics can be monotone, or have reliably computable long-term tendency, and can be closely approximated with the help of differential equations of some simple form. 

Also, modeling for a large number of dimensions of data becomes increasingly important. Our software assists in constructing relatively simple models of arbitrary dimensions. 

Constructing the differential equations from data helps us to obtain high precision (better than the traditional time series models) for the short-term prediction. Moreover, the dynamical system approach provides the global qualitative picture of the model's phase portrait, and shows the long-term tendencies of the process, underlying the time series realization. 

Our goal is to fit the data into the simple type of models with {\it trending flow} (discussed in~\cite{R4}, see definition~\ref{def-trending}). They allow us to use the phase portrait for the long-term predictions of the processes, and for better understanding of  the dynamics of the processes in nature, economics and social science. The phase portrait analysis of such models also helps us to understand various scenarios of the realization of the processes and associate them with the various trends of the phase portrait.

In many industries, users visit, join, or adopt a platform (such as content distribution service, payment system, or health insurance network) in order to access that platform's goods and services. 
One of the main factors of a platform's efficiency is the volume of users interacting through the platform. The dynamical system approach allows us to study the future behavior and the tendency of the trajectories of each group of platform users. It also helps to  increase the volume of users in the most cost efficient way.

We discuss these ideas with the help of the Wikipedia example. We use the data of Wikipedia's Readers, Contributors and Edits. 
In Section~\ref{section-error-analysis-2-dim} we show that the accuracy of the short-term prediction for the two-dimensional  dynamical system (with the Readers and Edits variables) is higher than the two-dimensional Auto-Regressive model. 
Then, (Section~\ref{3D-model}) we increase the dimensionality of the models, adding the Contributors data. This allows us to obtain a higher precision in one of the first two variables in VAR and in the dynamic model. Also, like in the two-dimensional case, the total error of the three-dimensional dynamical system is smaller than the total error of the three-dimensional VAR model.
In Section~\ref{sec-phase-portrate} we discuss the phase portraits of the two- and three-dimensional systems and the underlying characteristics of the process.

\section{Description of the dynamical system model and its short-term prediction}\label{section-error-analysis-2-dim}

In this section we discuss the questions of fitting models into the data and the accuracy of predictions. One of the important characteristics of a well-fitted ARIMA or a VAR model is randomness of their residuals, which can be measured with Box-Pierce test or a similar test. However, usually the ultimate goal is the forecast's accuracy.  Thus, in this work, our models' selection is such that it minimizes the error of the forecast, based on the knowledge of the past. 

The traditional time series models, such as VAR, ARIMA, Holt-Winters, are frequently used for short-term forecasts. We compare the precision of these forecasts to the dynamical system models' forecasts. The precision is measured as the sum of the squared distances between the true value taken from the testing data set and the value's short-term prediction, normalized by the magnitude of the true values. 

More specifically, we divide the data into two subsets, representing earlier time and later time. The first subset is used for the initial model construction. This model allows us to predict the first point of the second subset. We calculate the square of the distance between the predicted value and the true value. Then, we contribute this testing data point into the training set, re-evaluate the model, and predict the next value in the testing set. We calculate the square of the error in this second prediction and continue the process in the same way for the remaining testing points. Finally, we calculate the sum of the squared errors and divide it by the sum of the squares of the true values.

Comparing the errors of the models of various degrees, we choose the model that predicts the future values most accurately.
Using the Wikipedia data set, we found that the short-term prediction of the dynamical system (DS) model is more precise than the short-term prediction of the Vector Auto Regressive model (VAR).

What is the reason of the higher precision of the DS model? Both, the DS and VAR models, are designed to reconstruct the underlying ``true'' process, which generates the recorded time series. However, VAR (or a similar traditional) model attempts to smooth-out the noise, and it reconstructs the average behavior. The DS model uses the noise to its advantage, assuming that the noise does not affect the major law of the process, but shifts the realization of the process to some new trajectory (corresponding to a new initial condition) of the major process. Starting the next time-step prediction from the position that precisely corresponds to the current state, allows us to apply the law of the process to the true (non-averaged) current state. 

As an analogy, we can think of the Galileo Galilei experiment of dropping objects from the Pisa Tower during a strong storm. For the prediction of the height and velocity at the next time-step, we can apply the Free Fall law (the major process) to the current height and velocity (the new initial condition) at every time-step. Here, the Free Fall law can be estimated from the earlier data, and can be applied to the true current state. 

On the other hand, a traditional VAR (or similar) model can only use the current information for making more precise estimate of the average behavior.

Also, the precision of the DS prediction can be attributed to the fact that the model takes into account two conditions. The first one is the time series dependency of the later state on its earlier state  (derivatives of the model are calculated as the rate of change between two time-consecutive points). The second condition is the relation between the different coordinates of the multidimensional  state variable.

Different trajectories of the dynamical system may be associated with different external conditions that influence the system and force the transitions from one trajectory to another. 

In this work, the DS model is the system of differential equations, which have polynomial right-hand sides. 

First, we construct a DS model with the two variables: the number of Readers and the number of Edits. For the model construction and the error estimate we use the monthly Readers and Edits data of 2008-2019. 

We use the equations of the form~\eqref{seller-buyer-model}  (discussed in \cite{R3}) to model the traffic of Readers and Edits, $(x,y)$:
\begin{equation}\label{seller-buyer-model}
\left\{
\begin{array}{l}
x' = \epsilon_1 x +V_1(y),\\
y' = \epsilon_2 y +V_2(x).
\end{array}
\right. 
\end{equation}

Comparing the total errors of prediction for models having polynomial functions of degrees from 1 to 5, the equations with the degree 4 polynomial functions of the form~\eqref{Wiki-eqns} give the smallest error. We will call this model DS(4). 
\begin{equation}\label{Wiki-eqns}
DS(4):\ \left\{
\begin{array}{l}
x'=\epsilon_1 x +v_1y+v_2y^2+v_3y^3+v_4y^4,\\
y'=\epsilon_2 y +w_1 x- w_2x^2+w_3x^3 +w_4x^4.
\end{array}
\right.
\end{equation}

The best fitting autoregressive model is VAR(2). The total error comparison (table~\ref{error-table-pure-2D-data}) shows that DS(4) model gives more precise prediction than VAR(2) model.

\begin{equation}
\begin{tabular}{c|c|c|c}\label{error-table-pure-2D-data}
Model & Readers' predict. error & Edits' predict. error & Total error\\
\hline
DS(4)& .0051 & .0135 & .0186\\
VAR(2)& .0063 & .0134 & .0197
\end{tabular}
\end{equation}

The coefficients shown in~\eqref{Wiki-eqns-coefficients} were obtained for the model~\eqref{Wiki-eqns}, fitted in the entire data set. 
\begin{equation}\label{Wiki-eqns-coefficients}
\begin{array}{llllll}
\epsilon_1 & = & -0.3570, & \epsilon_2 &= & -0.2243\\
v_1 &= & -0.2637, & w_1 &= &1.2710\\
v_2 & = &6.9566, & w_2 &= & -6.9038\\
v_3 & = & -16.4522, & w_3 &= & 13.6668\\
v_4 & = & 11.0347, & w_4 &= & -8.6907.
\end{array}
\end{equation}

\begin{sidewaysfigure}
\centering
\includegraphics[width=\textwidth,height=\textheight,keepaspectratio]{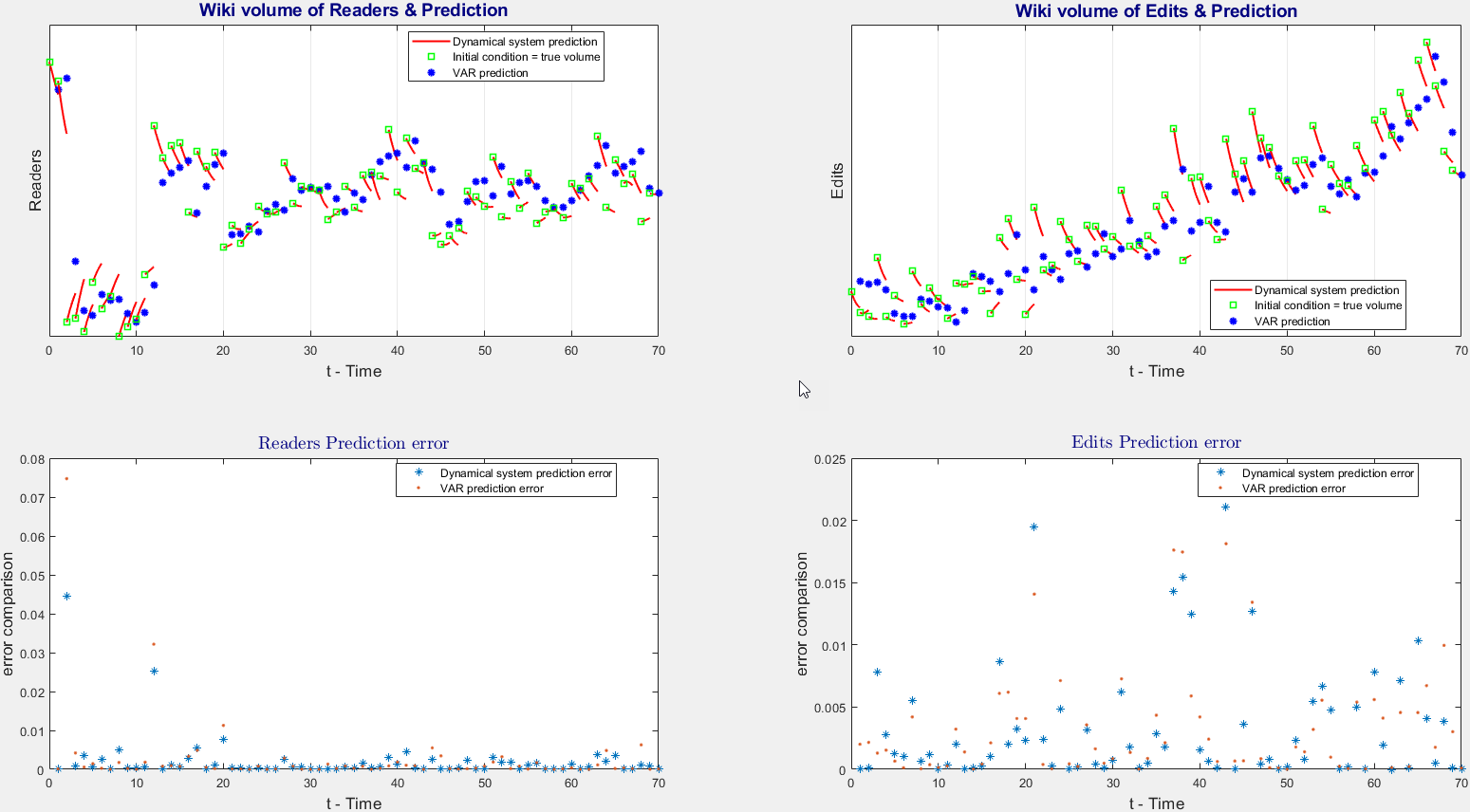}
\caption{The green dots in the upper Figures represent the number of Readers (left) and the number of Edits (right). The blue dots are the VAR(2) estimates of the Volumes. The red arrows start at the initial condition (current state, shown as the green dot) and end at the lag 1 predicted value. The lower Figures compare the errors of prediction via DS(4) and VAR(2) models for each points from the testing data set. The average error of the DS(4) is smaller than the one of the VAR(2) model.}
\label{single-variable-prediction-and-error}
\end{sidewaysfigure}

The Figure~\ref{single-variable-prediction-and-error} shows the predictions of each variable (Readers and Edits) with the help of Var(2) and DS(4). The figure also illustrates the size of the errors of the predictions. Each time-step prediction of the DS(4) starts at the current actual value of Edits and Readers and flows along the trajectory for the prediction at the next time. The total error of this prediction is smaller than the total error of the VAR(2) estimate.

We also consider a model, where the state variables are the Readers and Edits normalized by the growing (due to the growth of the number of Internet users in the world) potential number of Readers and the growing potential number of Edits. Factoring out the influence of the Internet development, creates a different dynamical model (with different global properties, discussed in Section~\ref{sec-phase-portrate}).

In this case the best fitting autoregressive model is also VAR(2) and the best dynamic model is also DS(4) model, having equations~\eqref{Wiki-eqns}. The total errors comparison (table~\ref{error-table-normalized-2D-data}) shows that DS(4) model gives more precise prediction than the VAR(2) model again.

\begin{equation}
\begin{tabular}{c|c|c|c}\label{error-table-normalized-2D-data}
Model & Readers' predict. error & Edits' predict. error & Total error\\
\hline
DS(4)&.0088 & .0116 & .0204\\
VAR(2)& .0109 & .0114 & .0223
\end{tabular}
\end{equation}

In this case, the DS(4) model's coefficients, based on the entire data set are shown in \eqref{Wiki-eqns-by-Internet-Users-coefficients}.

\begin{equation}\label{Wiki-eqns-by-Internet-Users-coefficients}
\begin{array}{llllll}
\epsilon_1 & = & -0.2677 & \epsilon_2 &= & -0.4655\\
v_1 &= & 2.3520 & w_1 &= & 1.1757\\
v_2 & = & -8.3986 & w_2 &= & -1.7697\\
v_3 & = &11.2901 & w_3 &= & 0.7948\\
v_4 & = & -5.1815 & w_4 &= & 0.0172
\end{array}
\end{equation}

\newpage
\vfill
\vfill
\clearpage

\section{The higher dimensionality of the dynamical system model}\label{3D-model}

Another advantage of the dynamical system's approach to the time series analysis is that the dynamical system's dimension can easily be modified. For example, for the Wikipedia traffic  estimate we may add one more variable: Contributors. The best fitting 3-dimensional dynamical system model is the DS(4), which has polynomials of degree 4 in its right side and negative $\epsilon_1=-0.4440,\ \epsilon_2=-0.3036,\ \epsilon_3=-0.3281$:

\begin{equation}\label{Wiki-3D-model}
DS(4):\ \left\{
\begin{array}{l}
x_1'=\epsilon_1 x_1 + P_4(x_2, x_3)\\
x_2'=\epsilon_2 x_2 +Q_4(x_1, x_3)\\
x_3' =\epsilon_3 x_3 +R_4(x_1, x_2)
\end{array}
\right.
\end{equation}

The origin is the fixed point, which has the eigenvalues $-1.7915 \pm 10.3549i$ and $2.5071$. 

As it is discussed in \cite{BK} and \cite{GJB}, the analysis of the large and complex time series data should benefit from the higher complexity models. The higher complexity can be associated with higher dimensionality. Sometimes additional variable can play a role, similar to the higher order of the model. Such new variable (new data) allows us to reduce the order of the model, or to achieve a higher precision.

The  total error of the 3-dimensional best fitting DS model (DS(4)) is smaller than the total error of best fitting VAR model VAR(3). Also, prediction of Edits in three-dimensional model is better than the one of the two-dimensional model (table~\ref{error-table-3D-data}):

\begin{equation}
\begin{tabular}{c|c|c|c|c}\label{error-table-3D-data}
Model & Readers' & Edits' & Contributors' & Total error\\
\hline
DS(4)& .0051&.0116 &.0020 &.0187\\
VAR(3)& .0053 &.0132 &.0031 & .0216
\end{tabular}
\end{equation}

\begin{sidewaysfigure}
\centering
\includegraphics[width=\textwidth,height=\textheight,keepaspectratio]{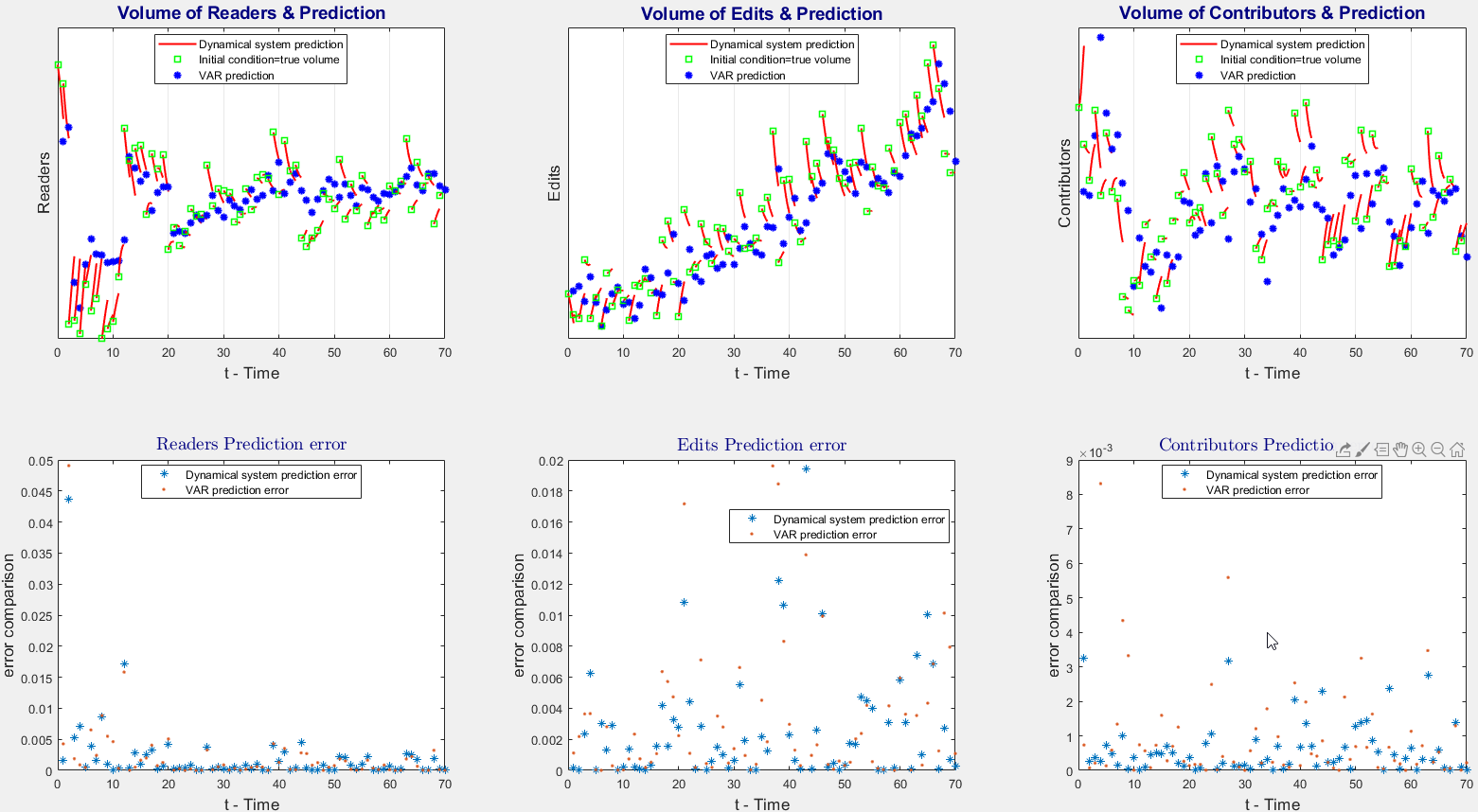}
\caption{The green dots in the upper Figures represent the volume of Readers, Edits and Contributors. The blue dots are the VAR(3) estimates of the Volumes. The red arrows start at the initial condition (current state, shown as the green dot) and end at the next time's value, predicted with DS(4). DS(4) makes a smaller error than VAR(3).}
\label{Wiki-error3D-DS4-VAR3}
\end{sidewaysfigure}

\vspace{.2in}
The Figure~\ref{Wiki-error3D-DS4-VAR3} shows the predictions of each variable (Readers, Edits and Contributors) with the help of Var(3) and DS(4). The figure also illustrates the size of the errors of the predictions. Each time-step prediction of the DS(4) starts at the current actual value of the variable and flows along the trajectory for the prediction at the next time. This total DS(4) error is smaller than the total VAR(3) error.

\newpage
\vfill
\vfill
\clearpage

\section{Phase portrait analysis}\label{sec-phase-portrate}

As it was discussed in Sections~\ref{section-error-analysis-2-dim} and \ref{3D-model} the dynamical system models provide an accurate short-term forecast. In this Section we will show that additionally, this approach allows us to see the main qualitative properties of the processes. The global picture of the flow provides information about trends, which can be associated with various scenarios of the process's realization.
 
Stationary points are important objects in the phase portrait analysis. They are the landmarks that organize the long-term behavior and describe the major characteristics of a process. The exact location of the stationary points is not informative; but the types of the stationary points explain the generic picture of the process and various scenarios of the process's realization. In our examples, we associate these scenarios with various basins of attraction of the fixed points. 

Phase portraits of some processes may have several basins of attraction. Fixed point in each basin shows tendency of the system in the long run. If the external conditions, regulations or incentives of the platform change significantly, the realization of the process may switch from one basin of attraction to another, and the two different fixed points can be associated with two different trends. Traditionally, this is modeled with the help of Intervention models, in which two different fixed points correspond to two different average behaviors. However, due to averaging, the Intervention models cannot carry as much information about the effects of the variation of the training data as the DS model. The DS models specifically reconstruct the effect of the rate of change of the training data.

Our software allows to fit time series data of any finite number of dimensions. We can also use a polynomial, trigonometric or square-root or other smooth function for the right side of the differential equations. 

However, we do not want to use high complexity equations, which generate a chaotic flow, because simpler equations fit the data sufficiently well and may not necessarily be improved via higher complexity. Also, nonchaotic nature of the flow permits more reliable analysis of the global properties and trends of the process. In \cite{R4} we defined the {\it trending flow}, which we use here for the Wikipedia traffic models.
\begin{definition}\label{def-trending}  
Assume that we are interested in the dynamics on the set $D$. Consider a system of differential equations defined on the domain $D$ (possibly, well-defined only on the interior of $D$). 
We will say that the system of differential equations has {\it trending flow on $D$} if its semiflow ($t\geq 0$) with any initial value in $D$ either converges to a fixed point in $D$, or escapes the domain $D$ (in finite time). 
\end{definition}

For the Internet platforms traffic, we use the system of differential equations (discussed in \cite{R4}) of the form~\eqref{seller-buyer-model-n-dim}: 
\begin{equation}\label{seller-buyer-model-n-dim}
\left\{
\begin{array}{l}
x_1 = \epsilon_1 x_1 +V_1(x_2,...,x_n),\\
x_2 = \epsilon_2 x_2 +V_2(x_1,x_3,...,x_n),\\
\vdots \\
x_n = \epsilon_n x_n +V_n(x_1,...,x_{n-1}),
\end{array}
\right. 
\end{equation}

where for all $i=1,..., n$, $x_i\geq 0$ and $V_i\geq 0$  on the domain of interest $D$. 

If $n=2$, the flow is trending. Also, if $\epsilon_i \geq 0$ ($i=1,.., n$), the flow is trending. See~\cite{R4}. There are examples of trending flow in the class of monotone systems that has been studied in depth by M.W. Hirsch. See, for example, \cite{H1, H2, H3, H4, H5, HS} and references therein. 

The dynamics of these processes can be viewed with the help of the software and analysis presented in this paper.

\subsection{Two-dimensional model for the Wikipedia's volume of Readers and Edits.}\label{section-phase-portrait-2D-pure}

In this subsection we discuss the planar dynamics defined by the variables Readers (estimated with the help of the 'number of page views' 2008-2019 data) and Edits (estimated with the helps of the `number of edits' 2008-2019 data). 

The best fitting parameters for this case are defined by the equations~\eqref{Wiki-eqns} with coefficients~\eqref{Wiki-eqns-coefficients} for the domain $D=[0,\infty)^2$. 

\begin{figure}
\centering
\includegraphics[scale=0.3]{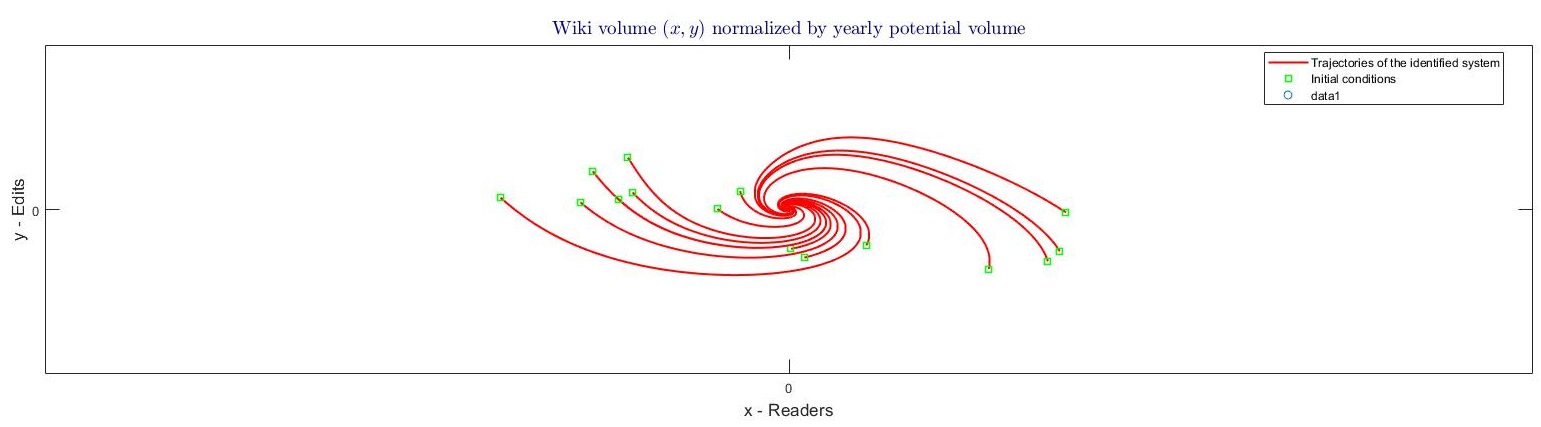}
\caption{In this figure, the flow belongs to the origin's basin of attraction. A small initial number of platform users oscillates and escapes the domain $D=[0, \infty)^2$ (Readers declines towards $0$). At the beginning of a platform formation, if the platform owners introduce some incentives, the flow may escape the the origin's basin and enter the basin of the fixed point $b$.}
\label{Wiki_phase_portrait_pure_data_neighbor0}
\end{figure}

In the domain of interest $D$, there are three positive fixed points: the origin, the positive (close to the origin) point $a=(a_1, a_2)$ and the positive fixed point $b=(b_1,b_2)$, both coordinates of which are significantly larger than the coordinates of the point $a$. The origin is a spiral attractor. The point $b$ is an attractor. The point $a$ is hyperbolic, through which the separatrix is passing. It separates (see Figure~\ref{Wiki_phase_portrait_pure_data_near_separatrix}) the origin's basin of attraction from the $b$ point's basin of attraction.
We can see that when the platform starts with a small positive number of users, its initial volume of users oscillates and escapes $D$ (Readers are vanishing). See Figure~\ref{Wiki_phase_portrait_pure_data_neighbor0}. The platform owners need to introduce some incentives at the beginning of platform's life. As soon as the volume of users reaches (via some jumps between trajectories, stimulated by externalities) the basin of the attracting stationary point $b$, it starts growing on its own. This model (based on the currently available data) suggests that no matter how big the volume of Internet users becomes, the volume of Wikipedia users remains bounded by the values of $b$. 

\begin{figure}
\centering
\includegraphics[scale=0.3]{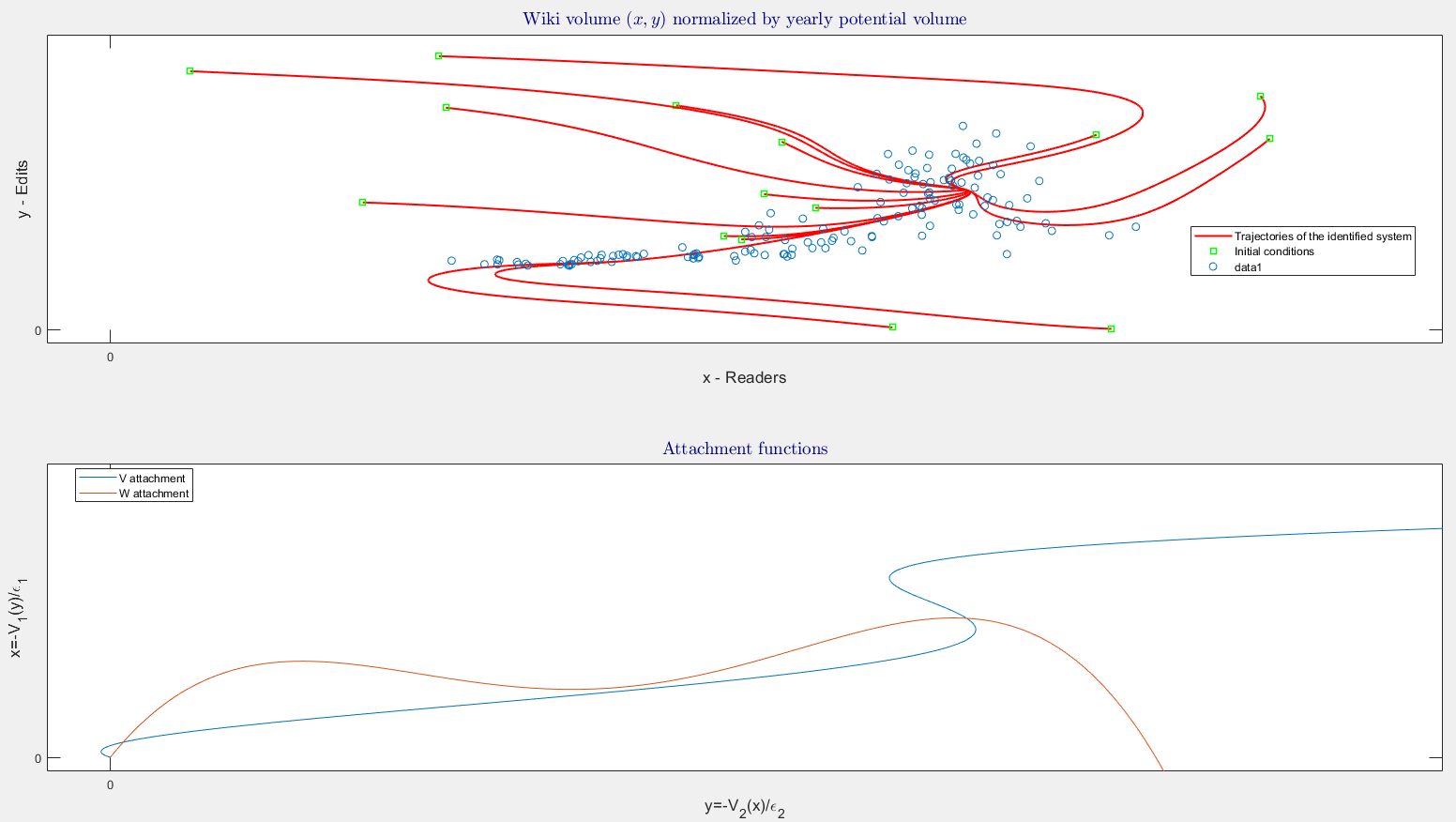}
\caption{The flow in the upper figure belongs to the basin of attraction of the point $b$.  If the volume gets above this fixed point (in one or both variables), it tends to eventually decrease towards the fixed point. The intercepts of the functions $-V(y)/\epsilon_1$ and $-V_2(x)/\epsilon_2$  (on the lower figure) define the location of the fixed points.}
\label{Wiki_phase_portrait_pure_data}
\end{figure}

\begin{figure}
\centering
\includegraphics[scale=0.3]{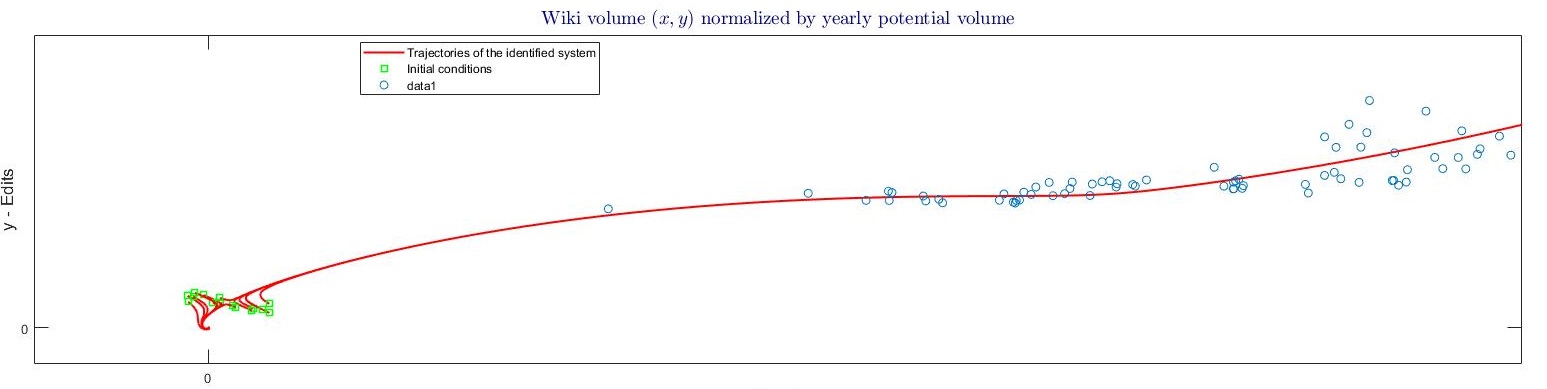}
\caption{The initial conditions (green squares) in this phase portrait are close to the separatrix (the division line between the two basins of attaction). When the volume of users is close to the separatrix, a small change may place it in either of the two basins: the positive sink (of the attractor shown in Figure~\ref{Wiki_phase_portrait_pure_data}) or in the basin of the origin's spiral sink shown in Figure~\ref{Wiki_phase_portrait_pure_data_neighbor0}. The initial conditions in this Figure demonstrate the two types of behavior: growth towards the point $b$ and decline towards $0$ Readers (escaping the domain $D$).}
\label{Wiki_phase_portrait_pure_data_near_separatrix}
\end{figure}

The function $V_1(y)$ is showing how the volume of Edits affects the traffic of Readers, and $V_2(x)$ is showing how the volume of Readers affects the volume of Edits. Clearly, $V_1$ and $V_2$ are positive valued functions on $D$ (see the bottom Figure~\ref{Wiki_phase_portrait_pure_data}), i.e. high number of Readers stimulates growth of Edits and high volume of Edits attracts more Readers. 

In this system of equations, $\epsilon_1 x$ is modeling the effect of the traffic of Readers on their rate of growth, and $\epsilon_2 y$ is modeling the effect of the traffic of Edits on their rate of growth. Both, $\epsilon_1,\ \epsilon_2\leq 0$. The sign of $\epsilon$ (the effect of the users of the same kind) was discussed in the papers~\cite{R1, R2, R3, R4}. The negative sign of this model can be explained by the ``edit wars'' (see~\cite{HKKR, HGMR}) on the Wikipedia platform. However, as discussed in the above referenced papers, this negative effect is small, if compared with the trading platforms, where sellers compete with each other for buyers, and buyers prefer low volume of buyers, which makes them more attractive to sellers and assures lower prices. This small negative effect creates small basin of attraction around the origin, and the platform owners need to provide small incentives when starting the platform, for the move into the basin of attraction of the positive fixed stationary point.

It can be shown that this two-dimensional model generates trending dynamics (for the details please see~\cite{R2}, \cite{R3}). The stationary points of the model~\eqref{seller-buyer-model} are located at the intercepts of the functions $x=-V_1(y)/\epsilon_1$ and $y=-V_2(x)/\epsilon_2$ (see the bottom of the Figure~\ref{Wiki_phase_portrait_pure_data}). The the flow shown in the phase portrait in the  Figure~\ref{Wiki_phase_portrait_pure_data} belongs to the basin of attraction of the point $b$. The basin of attraction of the origin and the neighborhood of the separatrix are shown in the Figure~\ref{Wiki_phase_portrait_pure_data_near_separatrix}. A more detailed picture of the behavior near the origin is shown in the Figure~\ref{Wiki_phase_portrait_pure_data_neighbor0}.

Thus, the phase portrait analysis suggests that there are two significant characteristics defined by the behavior (data) of this platform. These characteristics are associated with the fixed points of the dynamical system. The first one is the bounded growth (defined by the point $b$) of the platform's volume of users, independent of the growth of the Internet users. The second characteristic is the presence of the origin's basin of attraction. It demonstrates  that for a platform of this type, when the number of users is very low, some incentives may help to move the dynamics into the upper basin, and the platform's popularity will grow on its own.


\subsection{Factoring out the Internet growth effect}\label{section-phase-portrate-normalized-2D-data}

It is interesting to see what defines the main characteristics of dynamical system models discussed above, and how robust these characteristics are. Turns out that we can discuss the significance of the models' characteristics only when we compare them. If one of the characteristics becomes noticeably dominant, the signs of others may disappear (unless we increase the complexity of the model). 

\begin{figure}
\centering
\includegraphics[scale=0.3]{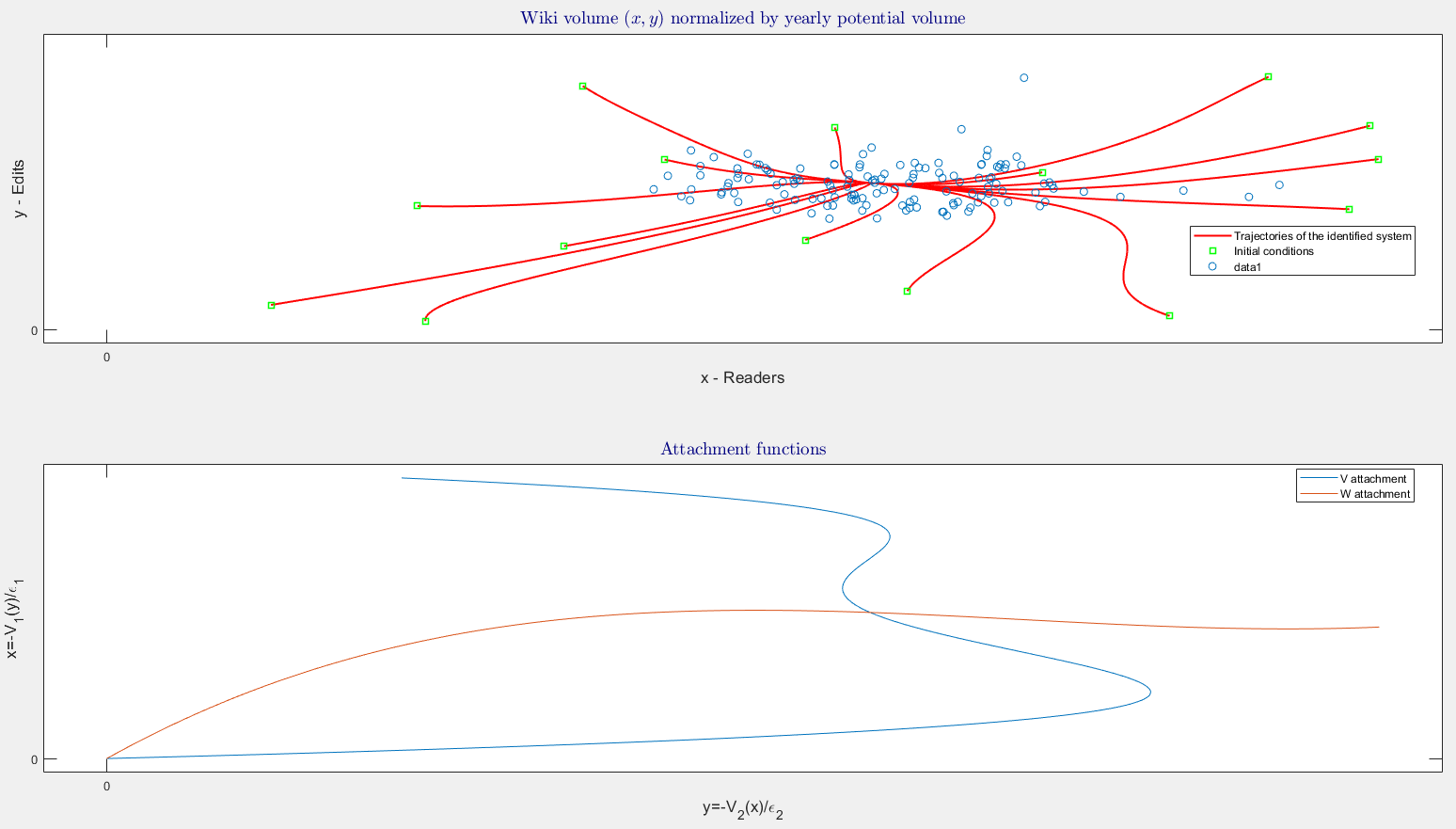}
\caption{We restrict our attention to the $[0,1]^2$ square. This square contains only 1 basin of attraction and 2 fixed points: $(0,0)$ and approximately $(.6, .5)$ (there are two more fixed points, located in the II and III quadrants -- outside of the domain of interest $D$). The intercepts of $-V(y)/\epsilon_1$ and $-V_2(x)/\epsilon_2$ show the location of the two fixed points in $[0,1]^2$.}
\label{phase-portrate-2D-normalized}
\end{figure}

In Chapter~\ref{section-error-analysis-2-dim} we discussed the fraction of all potentially available users. We assume that the potential number of Readers is defined by the number of Internet users, and the potential number of Edits additionally depends on the growth of R\&D. If we factor-out the influence of the Internet and R\&D growth, we can see that the most significant characteristic of the Wikipedia platform is the  bounded  growth of users, defined by the point, similar to the attractor $b=(b_1,b_2)$. The effect of the origin's basin of attraction becomes insignificant in comparison with the former characteristic. So, the model~\eqref{Wiki-eqns} with the coefficients~\eqref{Wiki-eqns-by-Internet-Users-coefficients} has only one basin of attraction, associated with the positive stationary point (see Figure~\ref{phase-portrate-2D-normalized}). The absence of the origin's basin of attraction may also be explained in this case by the ``small world'' idea. Namely, if there are no new Internet users, the permanent pool of the Internet users exchanges the information about the Wikipedia platform fast, and users motivate each other to join the platform. In this case, the ``edit wars'' are dominated by stronger effect of positive-valued $V_1$, $V_2$, the platform owners do not need to introduce incentives when they start the platform, and the platform grows on its own starting from any small initial volume of users. 

Thus, it is easier to start a new platform, if it is targeting a fixed (non-growing) pool of participants.

In this example of the 2-dimensional Wikipedia model (see equation~\eqref{Wiki-eqns} with coefficients~\eqref{Wiki-eqns-by-Internet-Users-coefficients}) on the domain $D=[0, 1]^2$, we consider the total number of Readers, divided by the total number of Internet users (UN data for the World); and the total number of Wikipedia Edits, divided by the total number of Internet users and adjusted by the rate of growth of the number of researches in the world (UN data).

\begin{figure}
\centering
\includegraphics[scale=0.3]{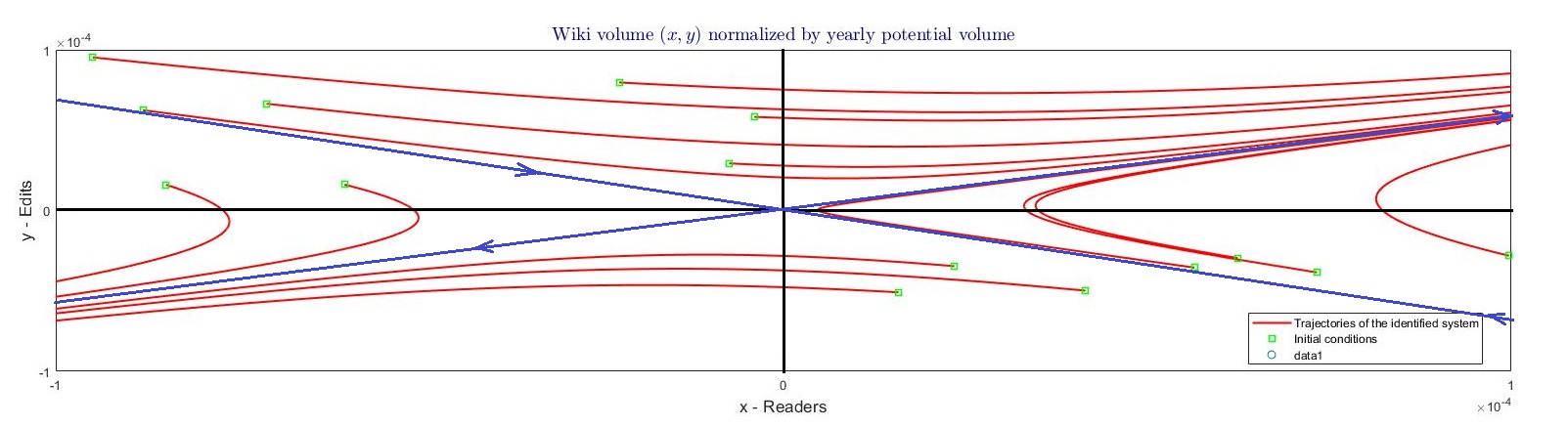}
\caption{Any small positive initial number of platform users tends to increase eventually, without incentives. I.e., all initial conditions of the first quadrant are driven up by the unstable manifold (shown in blue) in the first quadrant.}
\label{Wiki_phase-portrate-Poly-Order4-neighborhood0-edit}
\end{figure}

In order to understand the behavior of the flow of users, when the platform just starts, we analyze the equations near the origin. According to this model, the point $(0,0)$ is a hyperbolic fixed point (the eigenvalues of the linearized system are approximately $(-2.3159, 1.5620)$). The invariant (stable and unstable) manifolds of the linearized system are indicated as blue lines in the Figure~\ref{Wiki_phase-portrate-Poly-Order4-neighborhood0-edit}.  The  behavior of the volume of the platform users can be estimated as follows. Any small (non-zero)  initial volume of users will eventually increase in the direction of the unstable manifold. This is illustrated with the simulated initial conditions, which are very close to the origin\footnote{We show in the Figure~\ref{Wiki_phase-portrate-Poly-Order4-neighborhood0-edit} positive and negative initial conditions for a better understanding of the behavior near the origin, but we are interested only in the initial conditions in the first quadrant.}.  All trajectories tend to the non-zero fixed point. 

It is easy to see that this platform has trending flow.

This behavior differs Wikipedia from many popular platforms, which have strong negative same-side network effect (i.e., $\epsilon_1$ and $\epsilon_2$ are negative numbers of large magnitude and they are not dominated by $V_1$ and $V_2$). In the latter platforms, the origin is an attractor, and its basin of attraction has non-empty intersection with the first quadrant (which can be large). In this case (discussed in greater detail in~\cite{R5}), in order to escape this basin of attraction and to increase the volume of the platform's users, the platform owners need to introduce significant incentives, when they start their new platform. For example, a platform for interaction between sellers and buyers has a strong negative same-side network effect. It represents the competition between sellers (for buyers) and the competition between buyers for being the most desirable buyer (for sellers). Thus, if the same-side network effect is strong, it affects the flow near the origin (assuming that the origin is a stationary point) and makes the shrinking towards the origin to be a significant effect.

\subsection{Three-dimensional prediction for the Wikipedia platform.}\label{sec-phase-portrait}

In this section, we discuss the phase portrait of the tree-dimensional Wikipedia model with the new variable: Contributors (estimated with the help of the `number of contributors' data). We consider the system~\eqref{Wiki-3D-model} on the domain $D=[0,\infty)^3$. 

\begin{sidewaysfigure}
\centering
\includegraphics[width=\textwidth,height=\textheight,keepaspectratio]{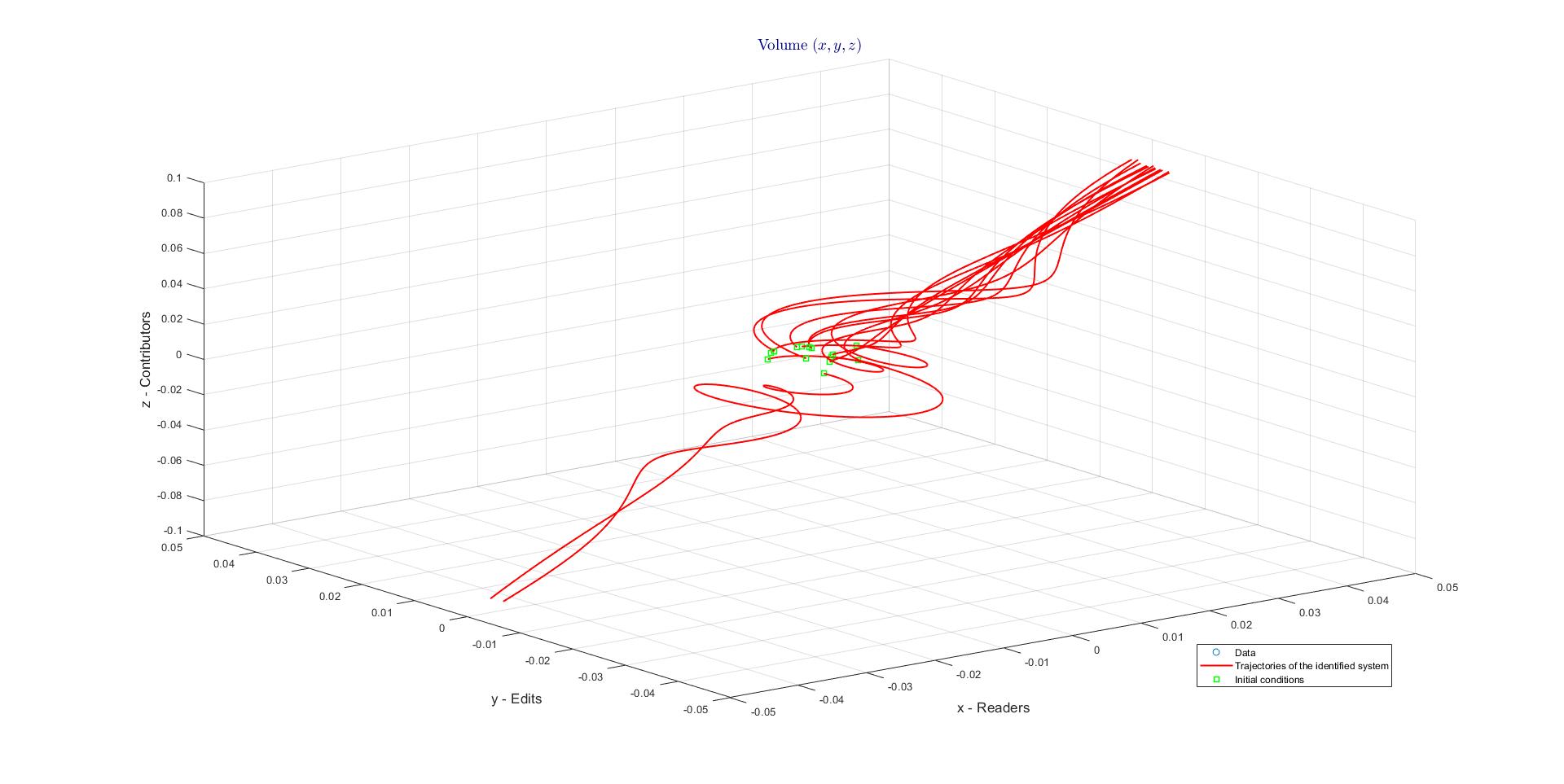}
\caption{The small initial volume of the 3-dimensional model.}
\label{Wiki3D-origin}
\end{sidewaysfigure}

In this phase portrait, the origin is the fixed point with one positive eigenvalue and with the conjugate pair of eigenvalues with the negative real part. So, near the origin, the trajectories spiral and either escape the domain $D$ or grow in the direction of the unstable manifold (associated with the increase of the number of Readers and Contributors). See the Figure~\ref{Wiki3D-origin}.

\begin{sidewaysfigure}
\centering
\includegraphics[width=\textwidth,height=\textheight,keepaspectratio]{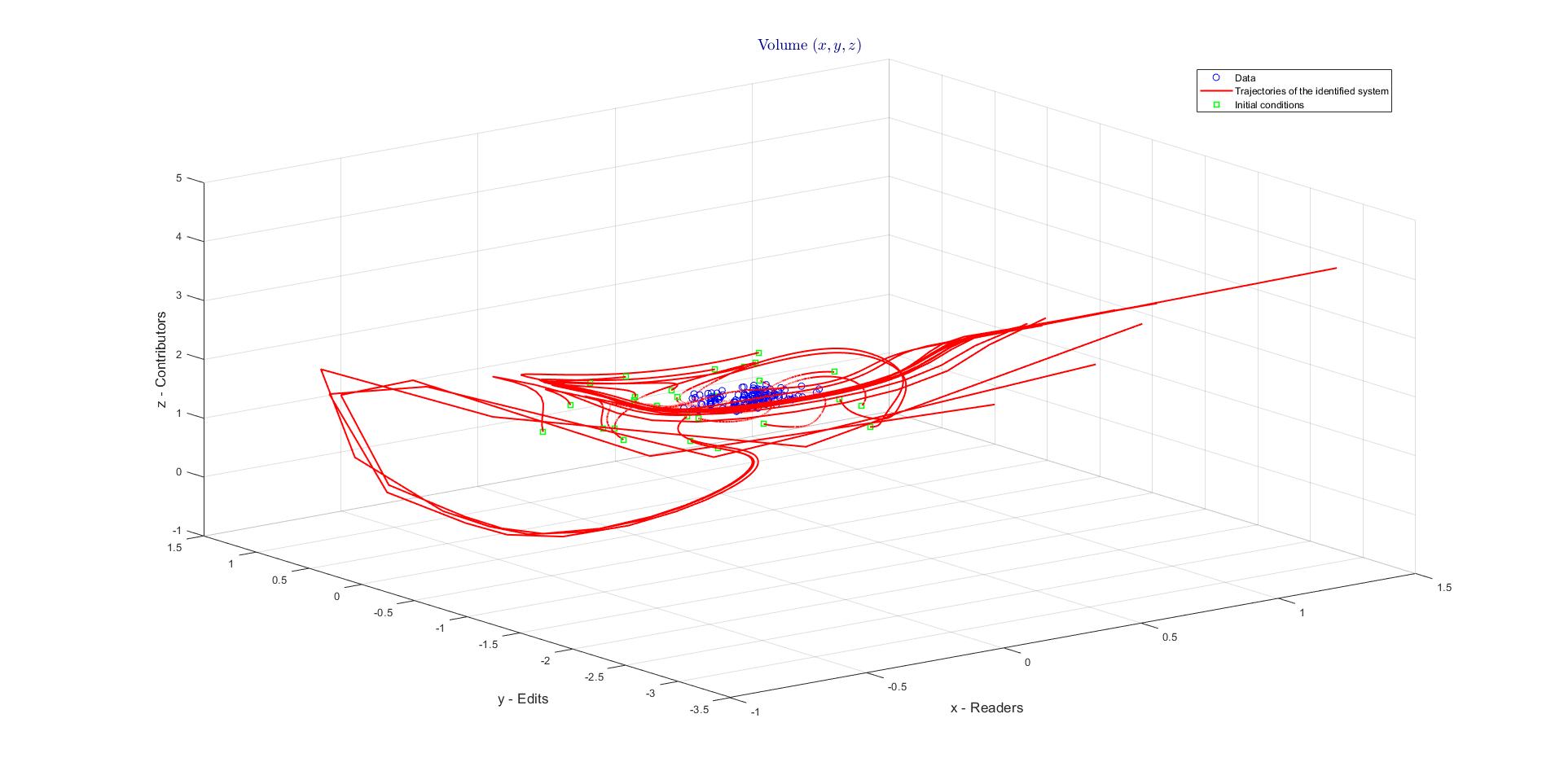}
\caption{The global view of the 3D-case. The tail in the right side of the picture shows the increase of the Readers and Contributors and the decline of the new Edits.}
\label{Wiki3D-global}
\end{sidewaysfigure}

In the interior of the first octant, this system has two more fixed points.  However, the behavior near the larger positive values (see Figure~\ref{Wiki3D-global}) is harder to analyze than the behavior near the origin.

There is one trend which we can see in the Figure~\ref{Wiki3D-global}: the volume of Readers and Contributors is growing, while the volume of Edits is decreasing. This can probably be explained by the fact that many of the subjects, in which Readers are interested, have already been contributed to the Wikipedia platform, and there is no demand for the new Edits. However, the Contributors will continue to correct some articles, and the Readers will be visiting the platform for the references.

It is interesting to investigate whether we can find a subset $D_1\subset D$, suitable for applications and such that we can prove that the dynamics on the $D_1$ subset is trending. 

It is also interesting to search for a simpler three-dimensional model, which would fit the data with the same or better accuracy, but such that we could prove that the model has trending dynamics that helps to explain various platform's development scenarios, associated with some significant changes in the external world, platform's policies and incentives.

\end{document}